\documentclass[preprint,12pt]{elsarticle}

\usepackage{amssymb}
\usepackage{amsmath}
\usepackage{lineno}
\usepackage{xcolor}
\usepackage{appendix}
\newproof{pf}{Proof}
\journal{Engineering Applications of Artificial Intelligence}

\begin{document}

\begin{frontmatter}


 \title{Cooperation dynamics between persuasion-prone and persuasion-averse populations under random sequential guidance}

\author[label1]{Lichen Wang}
\author[label1]{Shijia Hua}
\author[label1]{Yuyuan Liu}
\author[label1]{Siyu Liu}
\author{Linjie Liu\corref{cor1}\fnref{label1,label2}}
\ead{linjieliu1992@nwafu.edu.cn}
\cortext[cor1]{Corresponding Authors}

\affiliation[label1]{organization={College of Science, Northwest A\&F University},%Department and Organization
            city={Xianyang},
            postcode={712100}, 
            state={Shaanxi},
            country={China}}

\affiliation[label2]{organization={College of Economics \& Management, Northwest A\&F University},%Department and Organization
            city={Xianyang},
            postcode={712100}, 
            state={Shaanxi},
            country={China}}

\begin{abstract}
Social guidance, as a universal mechanism to promote collective cooperation, effectively maintains and enhances the level of population cooperation by correcting the behavior of defectors. However, existing guidance models have two key limitations: first, they assume that guidance actions are inevitably successful; second, they neglect information sharing among guides, leading to repeated guidance of the same successfully guided defector. Additionally, the psychological characteristics of individuals are often not fully considered in existing models. Here, we construct a stochastic sequence guidance model by introducing random ordering of guides and termination conditions for guidance actions to optimize the efficiency of guidance. Besides, we incorporate individual psychological traits such as persuasion-prone (where the success rate of persuasion increases with attempts) and persuasion-averse (where the success rate of persuasion decreases with attempts) into the modeling framework. We find that the introduction of random sequence guidance can not only effectively promote cooperation but also mitigate the second-order free-rider problem. Importantly, cooperation thrives with stronger persuasion, lower costs, and larger groups.

\end{abstract}




\end{frontmatter}

% \linenumbers

\section{Introduction}

Cooperation serves as a fundamental driver in the evolution of human societies. As ``super cooperators," humans have achieved evolutionary success primarily through their exceptional capacity for cooperation \cite{nowak2011supercooperators}. However, when the interests of individual agents conflict with group welfare, how cooperation emerges and stabilizes remains a central question in evolutionary biology and social sciences \cite{axelrod1981evolution,tanimoto2009promotion,perc2010coevolutionary,fu2017leveraging,perc2017statistical,wang2020eco}. The public goods game (PGG), as a commonly used framework for studying such dilemmas, demonstrates that under the assumption of rational self-interest, agents typically default to defection, leading to a collectively suboptimal outcome \cite{kurzban2001individual,szolnoki2010reward,tavoni2011inequality,charness2014identities,szolnoki2022tactical}. Evolutionary game theory provides important theoretical tools for understanding the emergence and maintenance of cooperative behavior within populations.

To address this dilemma, researchers have proposed various mechanisms to promote the emergence of cooperative behavior, primarily including reward, punishment, reputation, and exclusion methods \cite{brandt2003punishment,nikiforakis2008punishment,rand2011evolution,hauert2010replicator,sasaki2011replicator,sasaki2013evolution,sasaki2014rewards,liu2018evolutionary,yang2019reputation,quan2021comparison,wang2022replicator,quan2024reputation}. However, these measures all exhibit significant limitations: while incentive-based reward-punishment mechanisms can effectively promote cooperation, they incur high implementation costs and tend to induce second-order free-rider problems \cite{kiyonari2008cooperation}; reputation mechanisms rely heavily on sufficient information diffusion in social networks, demanding stringent requirements for information transparency \cite{kawakatsu2024mechanistic}; and although exclusion mechanisms can strongly suppress defection, their rigid exclusion criteria may inadvertently eliminate potential cooperators, ultimately undermining the collective payoff \cite{levitas1996concept,zhong2008cold}. Therefore, exploring more self-organized and highly inclusive cooperation mechanisms becomes particularly crucial.

The peer guidance approach works by having cooperative members actively persuade defectors to change their behavior, creating a self-sustaining system for maintaining cooperation \cite{groth2002commitment,watts2001career,leung2020mass}. This mechanism exhibits two distinctive advantages: first, its operation relies entirely on spontaneous interactions within the group without requiring intervention from external regulatory bodies; second, it maintains group cohesion through behavioral guidance rather than member exclusion. Taking community public projects as an example, actively participating members successfully convince reluctant neighbors to contribute through face-to-face communication, ultimately facilitating the realization of collective action.

However, existing models of guidance mechanisms commonly presume a static interaction framework wherein all guides continuously attempt to influence every defector until the end of a complete round \cite{10510635}. This oversimplification overlooks critical real-world dynamics: first, when a defector is successfully converted, the guidance activity should terminate immediately to conserve resources; secondly, the effectiveness of guidance can vary significantly due to social cognitive biases. For instance: certain populations might be persuasion-prone, with defectors becoming increasingly cooperative following numerous persuasions. Conversely, other populations might display persuasion-averse characteristics, wherein repeated persuasions may have counterproductive effects, diminishing the success rate of guidance. This work addresses these oversights by integrating these factors into an evolutionary game model and proposing a random sequence guidance model.

We adopt the replicator equation and Markov processes for theoretical and numerical analysis, respectively. The results show that the random sequential guidance mechanism can effectively promote cooperation in both infinite and finite populations. The parameter sensitivity analysis reveals that a higher initial guidance success rate is an important condition for achieving high-level cooperation. Meanwhile, the increase in group size, the enhancement of persuasion-prone tendencies or the reduction of persuasion-averse tendencies, as well as the decrease in guidance costs and information costs, all positively affect the final cooperation rate. Notably, we find that compared with populations that are persuasion-prone, populations that are persuasion-averse find it more difficult to achieve cooperation, and this difference is particularly evident when the initial guidance success rate is low.

\section{Related work}
Five core mechanisms have been established to explain the evolution of cooperation \cite{nowak2006five}. These mechanisms include direct reciprocity, indirect reciprocity, network reciprocity, group selection, and kin selection. Specifically, the direct reciprocity mechanism emphasizes repeated interactions between pairs of agents. When the expected number of interactions is sufficient, individuals will choose to cooperate in anticipation of future reciprocity \cite{hilbe2018partners}. Indirect reciprocity operates through a reputation system, through which cooperators can identify potential defectors based on interaction histories and thus adopt conditional cooperation strategies \cite{nowak2005evolution}. Network reciprocity theory shows that in structured populations, cooperators can form defensive clusters through spatial aggregation to effectively resist the invasion of defection strategies \cite{nowak1994spatial}. Group selection theory points out that when there is competition between groups, altruistic behavior may appear at the group level, but when intergroup competition disappears, cooperation may collapse within the group \cite{leigh2010group}. Kin selection explains the positive correlation between genetic relatedness and cooperative tendencies \cite{eberhard1975evolution}.

In governance mechanisms designed to promote cooperation and curb defection, punishing defectors or rewarding cooperators represent two fundamental strategies commonly employed by regulatory institutions. These strategies guide the evolution of group behavior toward cooperation by altering the payoff structure of games. Current research demonstrates that while both reward and punishment mechanisms can effectively enhance cooperation levels, each has inherent limitations: reward strategies typically require sustained high economic costs, whereas punishment mechanisms may trigger retaliatory or antisocial behaviors, thereby undermining group stability \cite{sigmund2001reward,hauert2007via,rand2009positive}. Existing studies indicate that adopting a ``First carrot, then stick" hybrid strategy (i.e., initially guiding cooperation through rewards, then implementing punishments against persistent defectors) can significantly reduce implementation costs while effectively achieving and maintaining cooperative equilibrium \cite{chen2015first}.

In addition, social exclusion can be viewed as a relatively severe punishment strategy \cite{li2015social}. This mechanism strictly removes defectors from the group, preventing them from free riding on public benefits, thereby effectively stopping free rider invasions and promoting cooperation. Existing research confirms that exclusion mechanisms implemented within a population (including both perfect exclusion and probabilistic exclusion) can significantly promote the emergence and maintenance of cooperative behavior. It is particularly noteworthy that when defectors appear in a population, exclusion mechanisms can simultaneously suppress second-order free riding problems. However, traditional exclusion models have a significant limitation: they typically assume that the decisions of excluders are completely independent, ignoring the exclusion outcomes achieved by other excluders. To address this, Li et al. \cite{li2016cooperation} proposed a random sequential exclusion model where each excluder implements exclusion with a fixed probability, and remaining exclusion behaviors stop immediately after any successful exclusion. Research has proven this dynamic exclusion mechanism performs excellently in maintaining cooperative stability \cite{li2016cooperation}. Although social exclusion mechanisms have been empirically demonstrated to effectively promote cooperation, their rigid implementation may produce dual adverse effects: not only undermining the group's social credibility but also potentially triggering retaliatory behaviors from excluded agents. Therefore, seeking mechanisms that provide opportunities for correcting defective behavior while maintaining team cooperation has emerged as a more optimal direction for institutional design.

As a flexible intervention method, the guidance mechanism has demonstrated unique evolutionary advantages in social practice. For instance, during global epidemic prevention efforts, altruistic agents consistently promote protective measures such as mask-wearing and social distancing. While these actions may momentarily decrease the individual advantages for the guides, they subsequently lead to a collective health improvement, thus mitigating infection risks and benefiting all, including the guides. Liu and Chen \cite{liu2024evolutionary} pioneered the integration of guidance strategies into the evolutionary game theory framework, devising two models: peer guidance and pool guidance. Based on the theoretical assumption that guidance behavior is bound to succeed (with a probability of 1), their studies not only validated the efficacy of guidance mechanisms in fostering cooperative behaviors but also highlighted its distinct effectiveness in curtailing second-order free rider issues \cite{liu2024evolutionary}.

Psychological studies have shown that the probability of behavioral change varies among agents as the frequency of guidance attempts increases. Specifically, persuasion-prone agents often resist changing their initial strategies during early stages, resulting in limited effectiveness of first attempts, but their cooperative tendencies increase significantly with cumulative guidance attempts. In contrast, persuasion-averse agents exhibit typical psychological reactance effects, demonstrating a negative correlation between persuasion frequency and compliance rates. This phenomenon is particularly evident in real-world scenarios such as consumer decision-making, where repeated persuasion attempts may trigger defensive resistance \cite{koch2017again}. Therefore, when exploring the impact of the guiding mechanism on the evolution of cooperative behavior, it is necessary to consider the termination conditions of the guidance as well as the psychological characteristics of the population.

The rest of this paper is organized as follows: Section III introduces the public goods game and establishes a guidance model. Section IV shows the comparative evolutionary results of replicator dynamics and stochastic dynamics. Section V summarizes the conclusions of this work and proposes some possible directions for future research.

\section{Models and methods}
\subsection{Public Goods Game}
We randomly select a group of $N$ individuals to participate in a public goods game. In this game, each individual must choose between two strategies: cooperation or defection. Those who choose to cooperate ($C$) must contribute a cost $c$ to the public pool, while those who choose to defect ($D$) do not contribute anything. The funds in the public pool are multiplied by a growth factor $r$ (where $1 < r < N$) and distributed equally among all agents. Although cooperation can generate higher total social benefits, rational agents tend to defect, constituting a social dilemma. In practice, introducing a guidance mechanism has become a common method for solving such dilemmas. The guidance provider can encourage cooperative behavior by persuading defectors.

\begin{figure*}[!t]
    \centering
    \includegraphics[width=1\textwidth]{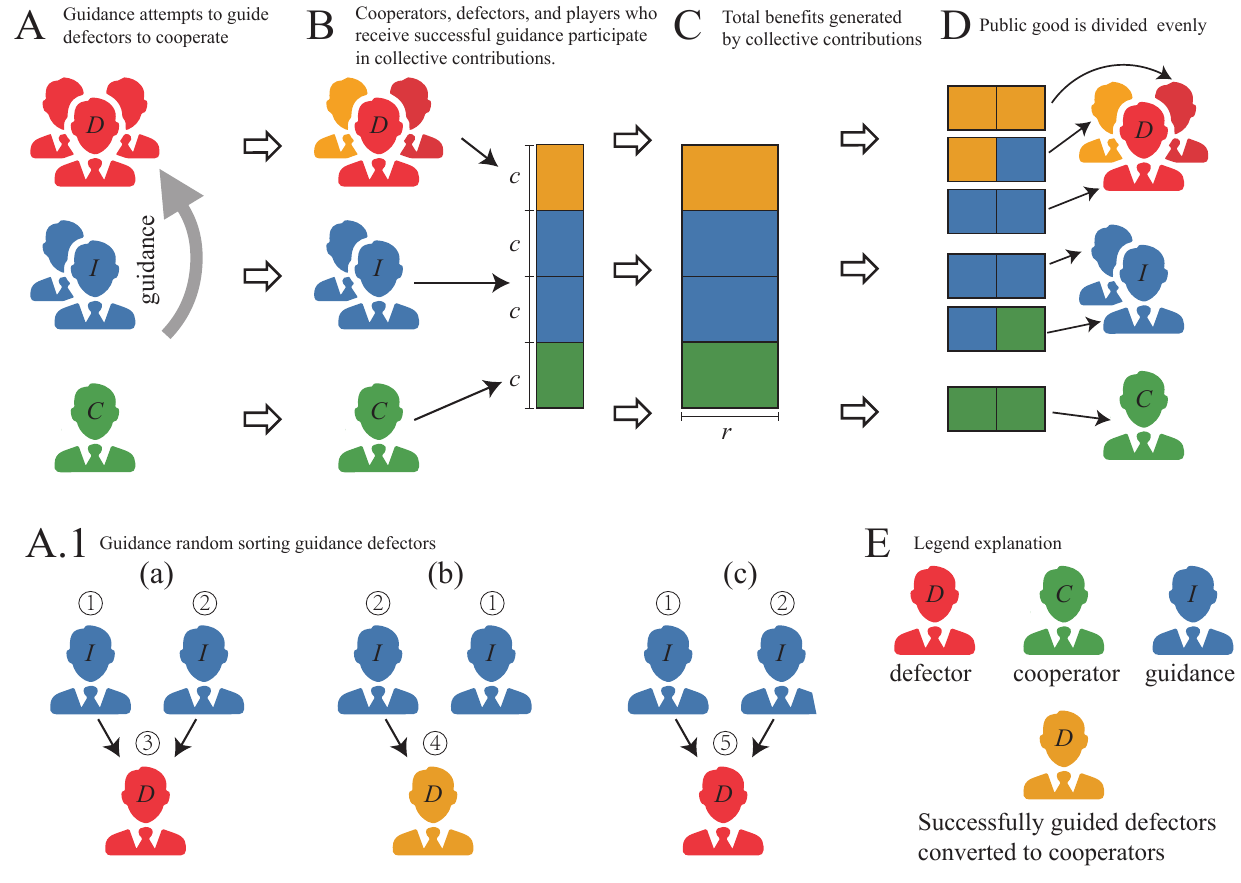}
    \caption{A schematic diagram illustrating how the guidance mechanism operates in public goods games. The interactions involve three defectors ($D$), two guides ($I$), and one cooperator ($C$). Panel A details the guidance process, where guides attempt to influence defectors to alter their strategies. This occurs in a randomized sequence, terminating either after the first successful guidance or once all guides have made an attempt. Specifically, sequences (a), (b), and (c) in Panel A.1 depict various orders of guidance and their outcomes; sequence A.1(a) and A.1(c) conclude when all guides have finished their attempts, whereas A.1(b) ends with the first successful guidance. Subsequently, in the donation phase (Panel B), all guides, cooperators, and successfully guided defectors contribute a fixed amount, $c$. These contributions are then pooled, yielding a return $r$ times the initial amount in Panel C. Finally, during the distribution phase (Panel D), the accumulated public goods are distributed equally among all group members, independent of their initial strategies.}\label{guidance map}
\end{figure*}

\subsection{Peer guidance}
In the framework of public goods games, each defector is regarded as a potential target for guidance. Assume that the guide's guidance order for the same defector is randomly sorted: when the target agent refuses to adopt a cooperative strategy in the guidance attempt, the subsequent guidance provider will continue to attempt guidance until the guidance is successful or all guidance providers have completed one guidance attempt (see Figure \ref{guidance map}). Empirical research in reality shows that the probability of successful guidance usually changes dynamically with the number of attempts. For this purpose, we set the initial probability of successful guidance as $\beta$, and let the change factor $\alpha$ describe the dynamic adjustment process of the probability. Specifically, the probability of successful guidance on the $k$-th attempt, $p_k$, is defined as
$$p_k = 
\begin{cases}
\alpha^{k-1}\beta, & \text{if } \alpha^{k-1}\beta \leq 1, \\
1, & \text{otherwise}.
\end{cases}$$
This model depicts two typical social behavior patterns: when $\alpha>1$, the probability of successful guidance increases with the number of attempts, indicating that agents tend to adopt cooperative strategies under multiple social influences. When $\alpha<1$, the probability of successful guidance decreases with the number of attempts, reflecting the resistance that agents may feel toward repeated persuasion.

For a randomly selected game group with $N_C$ cooperators, $N_D$ defectors, and $N_I=N-1-N_C-N_D$ guides (excluding the central individual), each guide can try to persuade the defectors to cooperate, and each guidance attempt incurs a fixed cost $\delta$. Once a defector is successfully guided, subsequent guidance attempts for that agent will be terminated to avoid wasting resources. For a single defector, the probability of failed guidance after a maximum of $N_I$ guidance attempts is:
$$\bar{p}(N_I) = \prod_{k=1}^{N_I}(1-p_k),$$
Correspondingly, the probability of successful guidance is
$$p(N_I) = 1 - \prod_{k=1}^{N_I}(1-p_k).$$

We assume that the sequence of guidance is randomized, with each guide's position $i$ in the sequence following a discrete uniform distribution $\mathcal{U}\{1, N_I\}$. With this configuration, the expected cost $C_R(N_I)$ incurred by a guide in persuading a single defector can be represented as follows:
\begin{equation*}
    C_{R}(N_I)= \sum_{i=1}^{N_I} \frac{1}{N_I}\prod_{h=1}^{i-1}(1-p_h) \delta,
\end{equation*}
where $\frac{1}{N_I}$ indicates the probability of a guidance provider being assigned to the $i$-th attempt. The term $\prod_{h=1}^{i-1}(1-p_h)$ denotes the probability that all previous $i-1$ attempts have failed. Notably, when $i=1$, this probability is set to 1.

The payoff of a strategy depends on the number of defectors successfully guided. If the number of successfully guided defectors is $N_W$, we can calculate the payoffs for individuals choosing to be cooperators, defectors, and guides, with the specific formulas as follows:
\begin{align*}
    \pi_C =
\begin{cases}
\dfrac{rc(N_C+1)}{N}-c, & \text{if } N_{I} = 0, \\
\dfrac{rc(N_C+1+N_{I}+N_W)}{N} - c, & \text{otherwise},
\end{cases}
\end{align*}
\begin{align*}
    \pi_D &=
    \begin{cases}
    \dfrac{rcN_C}{N}, & \text{if } N_{I} = 0, \\
    \begin{aligned}
    &p(N_I)\left(\frac{rc(N_C+N_{I}+N_{W}+1)}{N}-c\right)\\
     &\quad + \bar{p}(N_I)\frac{rc(N_C+N_{I}+N_{W})}{N},
    \end{aligned}
    & \text{otherwise},
    \end{cases}
\end{align*}

\begin{align*}
    \pi_{I} &=\dfrac{rc(N_C+N_{I}+N_W+1)}{N} - c - C_{R}(N_I+1) N_D - \theta,
\end{align*}
where $\theta$ quantifies the information cost required for guiding individuals to identify defectors within the population.

%In Figure \ref{guidance map}, we show the guidance process for a small group of size $N=6$ and the process of public goods donation and distribution.

\subsubsection{Replicator equations}
For an infinitely well-mixed population, let $x$, $y$, and $z$ represent the proportions of defectors, cooperators, and guides in the population, respectively. The dynamic evolution of these three strategies can be described by the following replicator equation:
\begin{equation}\label{system}
\begin{cases}
\dot{x}=x\left(f_{D}-\bar{f}\right), \\
\dot{y}=y\left(f_{C}-\bar{f}\right), \\
\dot{z}=z\left(f_{I}-\bar{f}\right),
\end{cases}
\end{equation}
where $f_D$, $f_C$, and $f_I$ represent the average fitness of defectors, cooperators, and guides, respectively. $\bar{f}=xf_D+yf_C+zf_I$ represents the average fitness of the population. The specific calculation is as follows:

\begin{align*}
     f_{C}=&\sum_{N_{D}=0}^{N-1}\sum_{N_{C}=0}^{N-N_{D}-1}\binom{N-1}{N_{D}}\binom{N-1-N_D}{N_{C}}x^{N_{D}}y^{N_C}z^{N_I}\\
    &\sum_{N_{W}=0}^{N_{D}}\binom{N_{D}}{N_{W}}p(N_I)^{N_{W}}\bar{p}(N_I)^{N_{D}-N_{W}}\pi_{C},
\end{align*}
\begin{align*}
     f_{D}=&\sum_{N_{D}=0}^{N-1}\sum_{N_{C}=0}^{N-N_{D}-1}\binom{N-1}{N_{D}}\binom{N-1-N_D}{N_{C}}x^{N_{D}}y^{N_C}z^{N_I}\\
    &\sum_{N_{W}=0}^{N_{D}}\binom{N_{D}}{N_{W}}p(N_I)^{N_{W}}\bar{p}(N_I)^{N_{D}-N_{W}}\pi_{D},
\end{align*}
\begin{align*}
     f_{I}=&\sum_{N_{D}=0}^{N-1}\sum_{N_{C}=0}^{N-N_{D}-1}\binom{N-1}{N_{D}}\binom{N-1-N_D}{N_{C}}x^{N_{D}}y^{N_C}z^{N_I}\\
    &\sum_{N_{W}=0}^{N_{D}}\binom{N_{D}}{N_{W}}p(N_I+1)^{N_{W}}\bar{p}(N_I+1)^{N_{D}-N_{W}}\pi_{I},
\end{align*}
where $N_I=N-1-N_C-N_D$.

In analyzing the frequency of cooperative behavior within a population, we focus on three principal elements: agents adopting cooperative strategies, those employing guidance strategies, and defectors who are successfully persuaded to adopt cooperative behaviors. At the equilibrium point $(x^\ast, y^\ast)$ of the system, the cooperation rate $\mathcal{C}^\ast$ and defection rate $\mathcal{D}^\ast$ are defined as follows:
\begin{align*}
    \mathcal{C}^\ast = &\sum_{N_{D}=0}^{N}\sum_{N_{C}=0}^{N-N_{D}}\binom{N}{N_{D}}\binom{N-N_D}{N_{C}}{x^\ast}^{N_{D}}{y^\ast}^{N_C}{z^\ast}^{N_I}\\
    &\times \begin{cases}
        \frac{N_C}{N}, &\text{if } N_{I} = 0, \\
        \sum_{N_{W}=0}^{N_{D}}\binom{N_{D}}{N_{W}}p(N_I)^{N_{W}}\bar{p}(N_I)^{N_{D}-N_{W}}\frac{N_C+N_I+N_w}{N},& \text{otherwise},
    \end{cases}\\
    \mathcal{D}^\ast = &\sum_{N_{D}=0}^{N}\sum_{N_{C}=0}^{N-N_{D}}\binom{N}{N_{D}}\binom{N-N_D}{N_{C}}{x^\ast}^{N_{D}}{y^\ast}^{N_C}{z^\ast}^{N_I}\\
    &\times \begin{cases}
        \frac{N_D}{N}, &\text{if } N_{I} = 0, \\
        \sum_{N_{W}=0}^{N_{D}}\binom{N_{D}}{N_{W}}p(N_I)^{N_{W}}\bar{p}(N_I)^{N_{D}-N_{W}}\frac{N_D-N_W}{N},& \text{otherwise},
    \end{cases}\\
    = &1-\mathcal{C}^\ast,
\end{align*}
where $z^\ast = 1-x^\ast-y^\ast$ and $N_I = N-N_C-N_D$.

\subsubsection{Markov process}
Consider a finite population of $Z$ agents, each with the option to adopt one of three strategies: cooperation ($C$), defection ($D$), or guidance ($I$). The system's state space $\mathcal{S}$ encompasses all potential combinations of these strategies: 
$$\mathcal{S}=\left\{\left(i_{C}, i_{D}\right) \in \mathbb{N}^{2} \mid 0 \leq i_{C}+i_{D} \leq Z\right\},$$ 
where $i_C$ and $i_D$ represent the number of cooperators and defectors, respectively, and the number of guides is determined by $i_I = Z - i_C - i_D$. Under the state $\mathbf{s} = (i_C, i_D)$, $N$ agents are randomly selected to interact. The fitness of cooperators, defectors, and guides are computed using the following hypergeometric distribution expressions, respectively:
\begin{align*}
    P_{C}=&\sum_{N_C=0}^{N-1}\sum_{N_{D}=0}^{N-N_{C}-1}\dfrac{\binom{i_C-1}{N_{C}}\binom{i_D}{N_{D}}\binom{Z-i_C-i_D}{N_I}}{\binom{Z-1}{N-1}}\\
    &\sum_{N_{W}=0}^{N_{D}}\binom{N_{D}}{N_{W}}p(N_I)^{N_{W}}\bar{p}(N_I)^{N_{D}-N_{W}}\pi_{C},
\end{align*}

\begin{align*}
    P_{D}=&\sum_{N_C=0}^{N-1}\sum_{N_{D}=0}^{N-N_{C}-1}\dfrac{\binom{i_C}{N_{C}}\binom{i_D-1}{N_{D}}\binom{Z-i_C-i_D}{N_I}}{\binom{Z-1}{N-1}}\\
    &\sum_{N_{W}=0}^{N_{D}}\binom{N_{D}}{N_{W}}p(N_I)^{N_{W}}\bar{p}(N_I)^{N_{D}-N_{W}}\pi_{D},
\end{align*}

\begin{align*}
    P_{I}=&\sum_{N_C=0}^{N-1}\sum_{N_{D}=0}^{N-N_{C}-1}\dfrac{\binom{i_C}{N_{C}}\binom{i_D}{N_{D}}\binom{Z-1-i_C-i_D}{N_I}}{\binom{Z-1}{N-1}}\\
    &\sum_{N_{W}=0}^{N_{D}}\binom{N_{D}}{N_{W}}p(N_I+1)^{N_{W}}\bar{p}(N_I+1)^{N_{D}-N_{W}}\pi_{I},
\end{align*}
where $N_I=N-1-N_C-N_D$.

At each time step, individuals adopting strategy $A$ may switch to strategy $B$ with a probability denoted as $p_{A \rightarrow B}$. This probability is defined using the Fermi function as follows \cite{szabo1998evolutionary}:
\begin{equation*}
	P_{(A \rightarrow B)}=\frac{1}{1+e^{\gamma(P_{A}-P_{B})}},
\end{equation*} 
where $P_A$ and $P_B$ describe the fitness of strategies $A$ and $B$, respectively, and $\gamma$ represents the strength of the influence of payoff differences on strategy selection. Considering the irrational behavior of individuals, we introduce a mutation probability $\mu$, where individuals randomly adjust their strategies with probability $\mu$. Under state $\mathbf{s}_i$, the one-step transition probability from strategy $A$ to strategy $B$ is given by:
\begin{equation*}
	T_{A \rightarrow B}^{\mathbf{s}_i}=(1-\mu)\left[\frac{i_{A}}{Z} \frac{i_{B}}{Z-1} \frac{1}{1+e^{\gamma\left(P_{A}-P_{B}\right)}}\right]+\mu \frac{i_{A}}{2 Z}.
\end{equation*}
Then, we can express the transition probabilities corresponding to all possible states that can be reached in one step from state $\left(i_C,i_D\right)$:
\begin{align*}
	&\left(i_{C}, i_{D}\right) \rightarrow\left(i_{C}+1, i_{D}\right)  = T_{I \rightarrow C}^{\mathbf{s}_i}, \\
	&\left(i_{C}, i_{D}\right) \rightarrow\left(i_{C}-1, i_{D}\right)  = T_{C \rightarrow I}^{\mathbf{s}_i}, \\
	&\left(i_{C}, i_{D}\right) \rightarrow\left(i_{C}, i_{D}+1\right)  = T_{I \rightarrow D}^{\mathbf{s}_i}, \\
	&\left(i_{C}, i_{D}\right) \rightarrow\left(i_{C}, i_{D}-1\right)  = T_{D \rightarrow I}^{\mathbf{s}_i}, \\
	&\left(i_{C}, i_{D}\right) \rightarrow\left(i_{C}+1, i_{D}-1\right)  = T_{D \rightarrow C}^{\mathbf{s}_i}, \\
	&\left(i_{C}, i_{D}\right) \rightarrow\left(i_{C}-1, i_{D}+1\right)  = T_{C \rightarrow D}^{\mathbf{s}_i},\\
	&\left(i_{C}, i_{D}\right) \rightarrow\left(i_{C}, i_{D}\right)  = 1-\sum_{U \neq U^{\prime}} T_{U \rightarrow U^{\prime}}^{\mathbf{s}_i},
\end{align*}
where $U, U' \in \{C, D, I\}$. Thus, we construct a strategy transition matrix of order $\frac{(Z+1)(Z+2)}{2}$. By solving for its normalized eigenvector with left eigenvalue 1, we obtain the stationary distribution of the Markov chain, denoted as $\mathbf{p}$.

Furthermore, the selection gradient at state $\mathbf{s}_i$ is defined as follows:
\begin{eqnarray}\label{Select gradient}
	\vec{\nabla}\left(\mathbf{\mathbf{s}_i}\right) & = & \left[\begin{array}{l}
		T_{I \rightarrow D}^{\mathbf{s}_i}+T_{C \rightarrow D}^{\mathbf{s}_i}-T_{D \rightarrow C}^{\mathbf{s}_i}-T_{D \rightarrow I}^{\mathbf{s}_i} \\
		T_{I \rightarrow C}^{\mathbf{s}_i}+T_{D \rightarrow C}^{\mathbf{s}_i}-T_{C \rightarrow I}^{\mathbf{s}_i}-T_{C \rightarrow D}^{\mathbf{s}_i}
	\end{array}\right].
\end{eqnarray}
This indicates the subsequent evolution path of the system under the current state $\mathbf{s}_i$.

After reaching a steady state, the cooperation rate $\mathcal{C}^\ast$ and defection rate $\mathcal{D}^\ast$ are quantified according to the following equations:
\begin{align*}\label{average adoption}
    \mathcal{C}^\ast &= C_{\mathbf{s}_i} \mathbf{p}_{\mathbf{s}_i},\\
    \mathcal{D}^\ast &= D_{\mathbf{s}_i} \mathbf{p}_{\mathbf{s}_i}
    = 1-\mathcal{C}^\ast,
\end{align*}
where $\mathbf{p}_{\mathbf{s}_i}$ is the probability of staying in state $\mathbf{s}_i$. $ C_{\mathbf{s}_i}$ and $ D_{\mathbf{s}_i}$ represent the proportions of cooperation and defection in state $\mathbf{s}_i$, respectively. The specific form is as follows:
\begin{align*}
       C_{\mathbf{s}_i} = &\sum_{N_C=0}^{N}\sum_{N_{D}=0}^{N-N_{C}}\dfrac{\binom{i_C}{N_{C}}\binom{i_D}{N_{D}}\binom{Z-i_C-i_D}{N_I}}{\binom{Z}{N}}\\
       &\times \begin{cases}
        \frac{N_C}{N}, &\text{if } N_{I} = 0, \\
        \sum_{N_{W}=0}^{N_{D}}\binom{N_{D}}{N_{W}}p(N_I)^{N_{W}}\bar{p}(N_I)^{N_{D}-N_{W}}\frac{N_C+N_I+N_w}{N},& \text{otherwise},
        \end{cases}\\
       D_{\mathbf{s}_i} = &\sum_{N_C=0}^{N}\sum_{N_{D}=0}^{N-N_{C}}\dfrac{\binom{i_C}{N_{C}}\binom{i_D}{N_{D}}\binom{Z-i_C-i_D}{N_I}}{\binom{Z}{N}}\\
       &\times \begin{cases}
        \frac{N_D}{N}, &\text{if } N_{I} = 0, \\
        \sum_{N_{W}=0}^{N_{D}}\binom{N_{D}}{N_{W}}p(N_I)^{N_{W}}\bar{p}(N_I)^{N_{D}-N_{W}}\frac{N_D-N_W}{N},& \text{otherwise},
    \end{cases}
\end{align*}
where $N_I=N-N_C-N_D$.

Complete theoretical analysis of two-strategy evolutionary dynamics, including both deterministic and stochastic formulations, is provided in the Supplementary Information.

\section{Results}
\begin{figure*}[!t]
    \centering
    \includegraphics[width=1\textwidth]{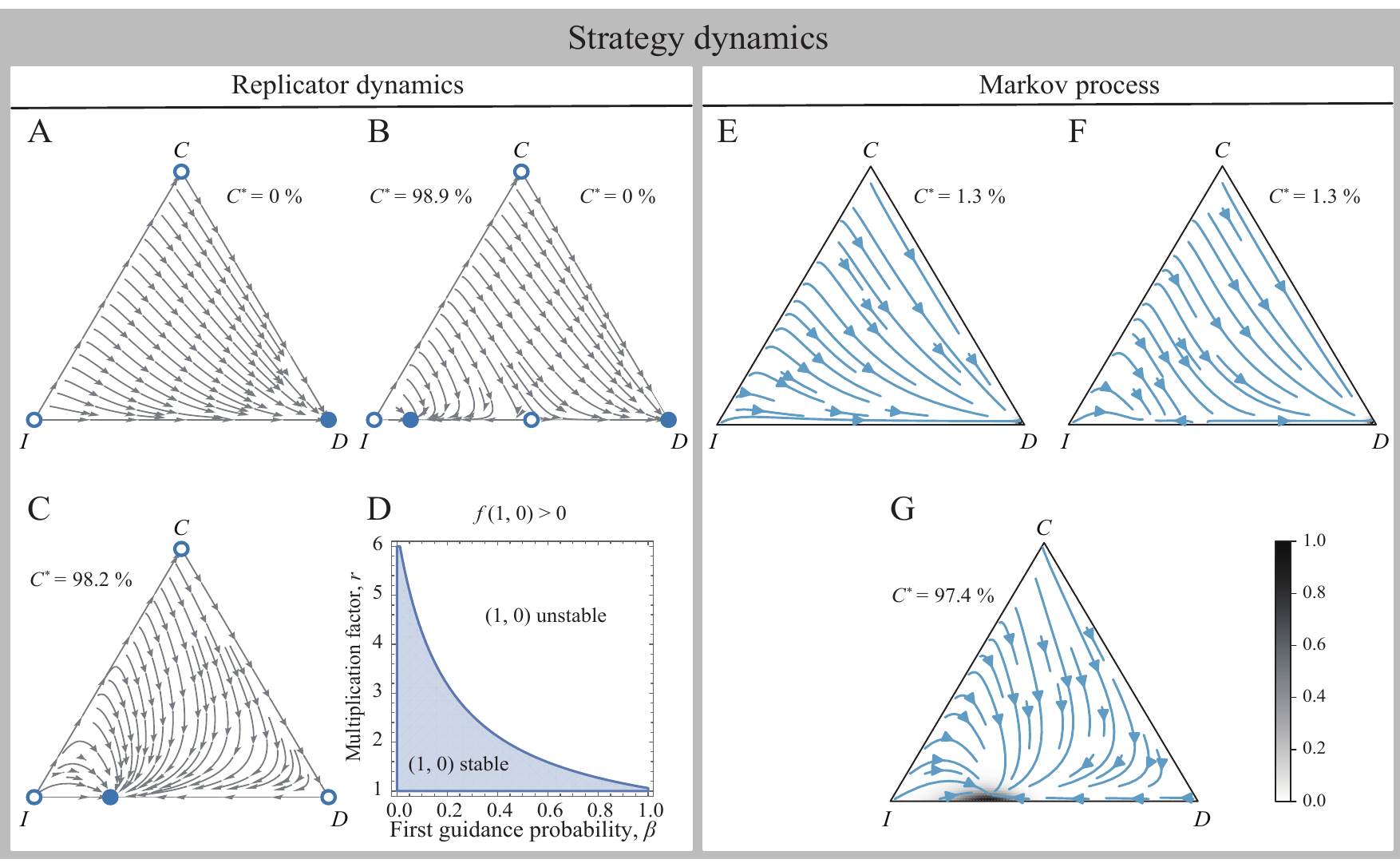}
    \caption{The evolutionary dynamics are characterized by both replicator equations and Markov processes. Panels A to C show the phase diagram analysis from replicator dynamics: under conditions of low initial guidance probability ($\beta=0.1$) and weak persuasion propensity ($\alpha=1.6$), the system converges to a full defection equilibrium (Panel A); when maintaining $\beta=0.1$ while increasing persuasion propensity to $\alpha=2$, the system exhibits bistable characteristics, where raising the initial proportion of guides significantly enhances cooperation, achieving a steady-state cooperation rate of $\mathcal{C}^\ast = 98.9$ (Panel B); whereas with high initial guidance probability ($\beta=0.5$) but relatively low persuasion propensity ($\alpha=1.1$), the system can still maintain high-level stable cooperation at $98.2$ (Panel C). Panel D further reveals a negative correlation between social dilemma intensity ($r$-value) and the required initial guidance success rate to avoid complete defection. Panels E to G display the evolutionary trajectories from the stochastic dynamics, where color intensity represents the stationary distribution across system states. The parameters in panels A, B, and C are $N=6$, $c=1$, $r=3$, $\delta=0.01$, and $\theta=0.01$. In panel D, parameters are $N=6$, $c=1$, $\delta=0.01$, and $\theta=0.01$. For panels E, F, and G, the values of $\beta$ and $\alpha$ are identical to those in panels A, B, and C, respectively. The additional parameters are $Z=100$, $N=6$, $c=1$, $r=3$, $\delta=0.01$, $\mu=1/Z$, $\gamma=10$, and $\theta=0.01$.}\label{Strategy dynamics}
\end{figure*}

Now we focus on exploring the determinacy and stochastic dynamics of systems in both infinite and finite well-mixed populations. The replicator dynamics demonstrate three characteristic outcomes: defector dominance, bistability (between defector dominance and defector-guide coexistence), and defector-guide coexistence. As shown on the left side of Figure \ref{Strategy dynamics}, there are three corner equilibrium points, namely, $(0,0)$, $(0,1)$, and $(1,0)$, where the first two are unstable and the last one is stable when $\theta+\delta(N-1)+c-\frac{rc[(N-1)\beta+1]}{N}>0$. It indicates a scenario where the cooperation rate ultimately reaches zero and all agents opt for defection (see Figure \ref{Strategy dynamics}A). In the Supplementary Information, we prove that, except for the equilibrium point $(1,0)$ corresponding to complete defection, stable equilibrium points can only exist at the boundary where defection and guidance coexist. We denote it as $(x^\ast,0)$. Figure \ref{Strategy dynamics} B illustrates a bistable situation where, besides the stable equilibrium of full defection, a second stable equilibrium exists comprising guides and defectors, which significantly enhances cooperation levels within the population. Figure \ref{Strategy dynamics} C depicts a scenario in which only guidance and defection strategies coexist at equilibrium, achieving a high cooperation rate. This is attributed to the high probability $\beta$ of successful initial guidance by the guidance player, which facilitates the replacement of the complete defection strategy by the guidance strategy, thus rendering the equilibrium point $(1,0)$ unstable. Figure \ref{Strategy dynamics} D examines the effects of social dilemma intensity and the probability of initial guidance on the equilibrium point $(1,0)$. It shows that as the social dilemma intensity decreases (i.e., $r$ increases), the probability of initial guidance required to prevent a full defection equilibrium decreases.

The right side of Figure \ref{Strategy dynamics} shows the stochastic dynamics of the system. When the initial guidance probability is low coupled with weak persuasiveness, the population evolves toward complete defection dominance (see Figure \ref{Strategy dynamics} E). Although Figure \ref{Strategy dynamics} F shows that when there are enough guides, the system may evolve to a state where guidance and defection coexist, the system is still more likely to remain in a state of high defection. Nevertheless, under elevated initial guidance probability ($\beta = 0.5$) despite moderate persuasion tendency ($\alpha = 1.1$), the system can achieve stable cooperation maintenance at high levels. These evolutionary dynamics are similar to the results revealed by the replicator equation.

\begin{figure*}[!t]
    \centering
    \includegraphics[width=1\textwidth]{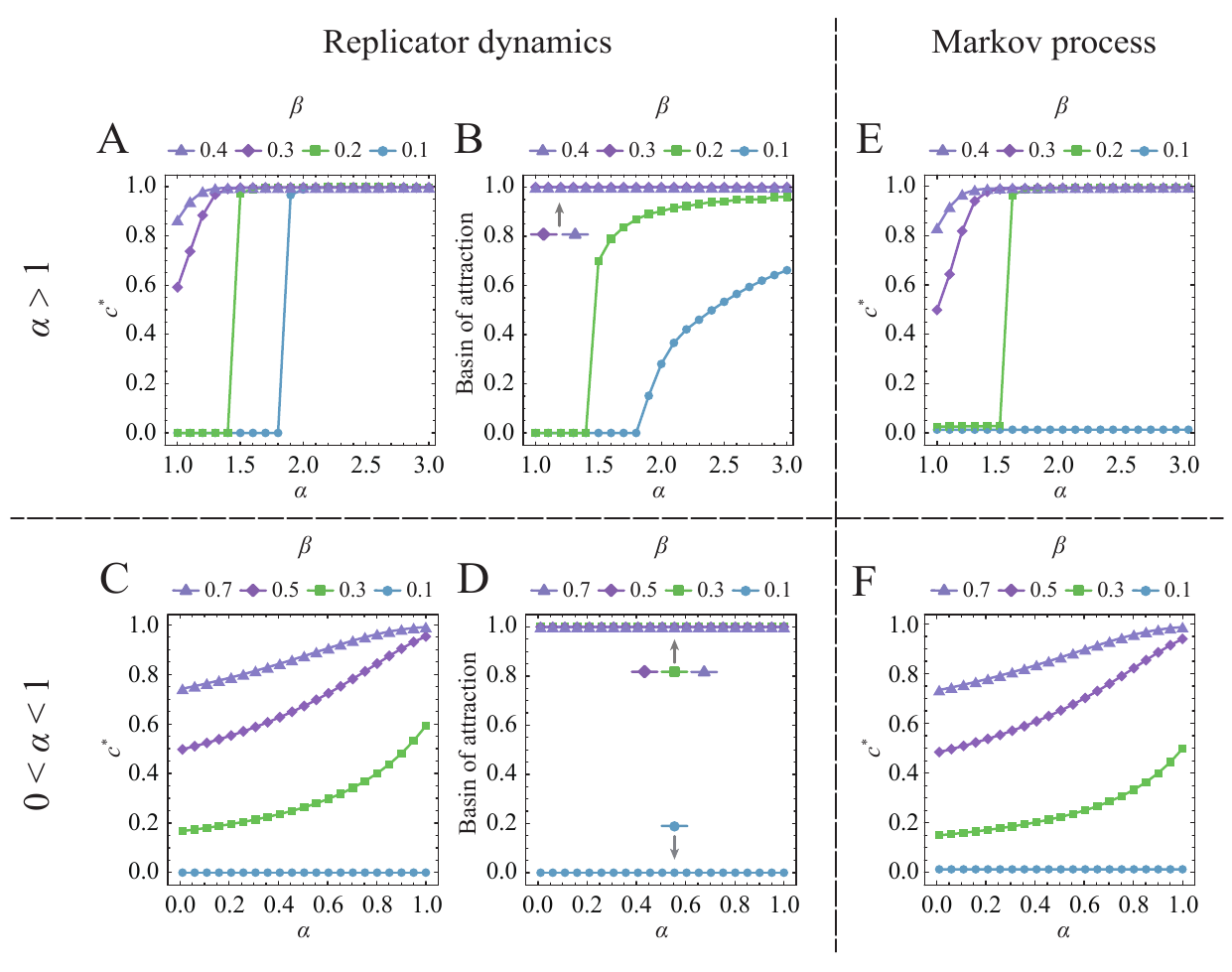}
    \caption{The effects of parameter $\alpha$ under varying initial guidance success rates $\beta$ on two characteristic population types: persuasion-prone populations ($\alpha>1$) and persuasion-averse populations ($0<\alpha<1$). Panels A and C display the cooperation rates $\mathcal{C}^\ast$ at the stable equilibrium point $(x^\ast,0)$ for different population types. Panels B and D provide quantitative analysis of the basin of attraction properties for this equilibrium. Panels E and F illustrate how parameter $\alpha$ modulates stochastic dynamics. The integrated results demonstrate that increasing the values of parameters $\alpha$ and $\beta$ significantly enhances cooperation levels. The parameters for panels A–D are $N=6$, $c=1$, $r=3$, $\delta=0.01$, and $\theta=0.01$. The parameters for panels E–F are $Z=100$, $N=6$, $c=1$, $r=3$, $\delta=0.01$, $\mu=1/Z$, $\gamma=10$, and $\theta=0.01$.}\label{alpha}
\end{figure*}

Population behavior often exhibits one of two distinct characteristics: persuasion-prone or persuasion-averse. These tendencies profoundly influence the evolutionary dynamics of the system. As depicted in Figure \ref{alpha}, we employ modeling frameworks, specifically replicator dynamics and Markov processes, to analyze the impact of parameter $\alpha$ and the probability of initial guidance success $\beta$ on population cooperation levels and the basin of attraction at the stable equilibrium point $(x^\ast,0)$. A basin of attraction value of zero suggests that the equilibrium point $(x^\ast,0)$ is unstable. 

When the parameter $\alpha>1$, individuals show clear persuasion-prone characteristics, meaning defectors are more likely to be persuaded to participate in cooperation as guidance attempts increase. In replicator dynamics analysis, when $\beta$ takes small values and $\alpha$ is also small, the presence of many guides is insufficient to significantly increase the probability of successful guidance. Consequently, defector strategies dominate, driving the system toward a ``tragedy of the commons" characterized by full defection. Conversely, when $\alpha$ is high, increasing the number of guides raises the probability of successful persuasion, promoting the emergence of a stable equilibrium $(x^\ast,0)$ where guides and defectors coexist (Figure \ref{alpha} A). Moreover, as $\alpha$ increases, the basin of attraction of this equilibrium also expands, indicating guides play a key role in promoting high cooperation levels. Notably, when $\beta$ is large, $(x^\ast,0)$ becomes the only stable equilibrium (Figure \ref{alpha} B). The results in finite populations agree with replicator dynamics, confirming that larger $\beta$ and $\alpha$ values remain key factors for maintaining high cooperation levels (Figure \ref{alpha} E). When $\beta$ is not too small, the stochastic dynamics match replicator dynamics, but when $\beta$ is small, even with large $\alpha$, the stochastic dynamics show low cooperation rates in the population. Our results demonstrate that larger values of $\beta$ and $\alpha$ are crucial to maintain higher level of cooperation.

When $0<\alpha<1$, agents exhibit distinct repeated-persuasion aversion characteristics, and this aversive effect intensifies as $\alpha$ decreases. The results demonstrate that for populations with persuasion-averse tendencies, the effectiveness of initial guidance is particularly crucial. Specifically, when the initial guidance success rate ($\beta$) remains high, moderately increasing $\alpha$ can significantly enhance the level of cooperation (Figure \ref{alpha} C and F); conversely, when $\beta$ is relatively low, even with weak persuasion aversion, the inadequate initial guidance will still lead to systemic collapse of cooperative behavior (Figure \ref{alpha} C, D, and F).

\begin{figure*}[!t]
    \centering
    \includegraphics[width=1\textwidth]{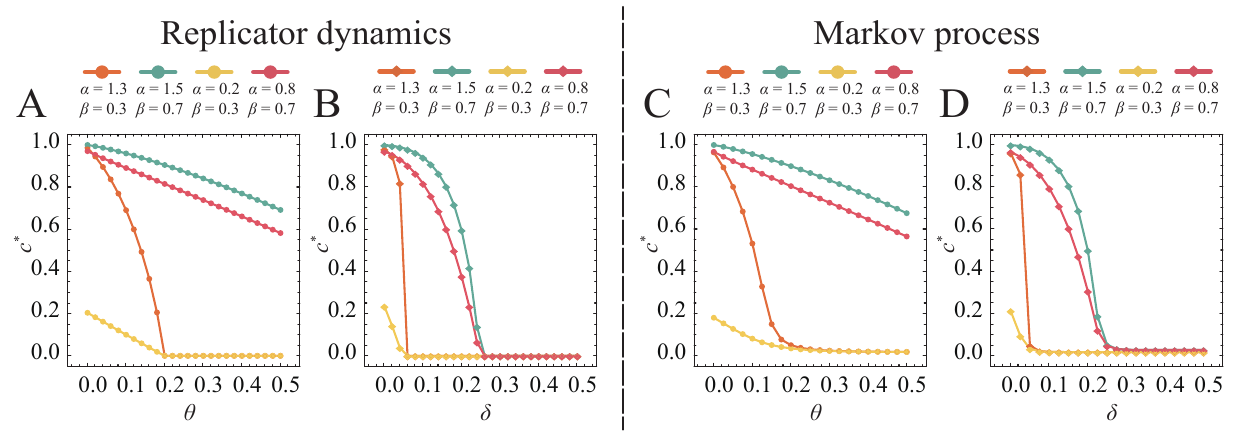}
    \caption{The effects of information costs $\theta$ and guidance costs $\delta$ on the cooperation rate. Two typical combinations of scenario parameters are set up: (1) poor guidance environment conditions ($\alpha=1.3$, $\beta=0.3$; $\alpha=0.2$, $\beta=0.3$); (2) favorable guidance environment conditions ($\alpha=1.5$, $\beta=0.7$; $\alpha=0.8$, $\beta=0.7$). The results show that both an increase in information acquisition cost $\theta$ and an increase in guidance cost $\delta$ will lead to a decline in the level of cooperation. The remaining parameters for panel A are $N=6$, $c=1$, $r=3$, and $\delta=0.01$; $N=6$, $c=1$, $r=3$, and $\theta=0.01$ in panel B; $Z=100$, $N=6$, $c=1$, $r=3$, $\delta=0.01$, $\mu=1/Z$, and $\gamma=10$ in panel C; $Z=100$, $N=6$, $c=1$, $r=3$, $\theta=0.01$, $\mu=1/Z$, and $\gamma=10$ in panel D.}\label{theta_delta}
\end{figure*}

Unlike traditional cooperation and defection strategies, the implementation of guidance strategies depends on the accurate identification of defectors within a group. However, due to the possibility of defectors engaging in disguise behavior, the cost of information identification $\theta$ is often very high, which may weaken the evolutionary advantage of guidance strategies, thereby affecting the overall level of cooperation within the population. As shown in Figure \ref{theta_delta}A and Figure \ref{theta_delta}C, increasing information costs significantly weaken the competitive advantage of guidance strategies, ultimately leading to a decline in the system's cooperation rate. Notably, under poor guidance conditions, when information costs exceed a critical threshold, the system inevitably falls into a ``tragedy of the commons" state where all individuals tend to adopt defection strategies, resulting in complete collapse of cooperation.

Guidance cost, as a key constraint that affects individuals' adoption of guidance strategies, specifically refers to the resource expenditure (including time and monetary costs) that guides must invest in order to reverse the strategic choices of defectors. As shown in Figure \ref{theta_delta}B and D, guidance costs have a decisive impact on cooperation levels: when costs exceed a critical threshold, full defection becomes the system's sole evolutionarily stable state. Notably, in favorable guidance environments (e.g., with higher $\alpha$ and $\beta$ values), individuals exhibit greater willingness to bear costs and are more inclined to invest substantial resources to guide defectors toward strategy change.
\begin{figure*}[!t]
    \centering
    \includegraphics[width=1\textwidth]{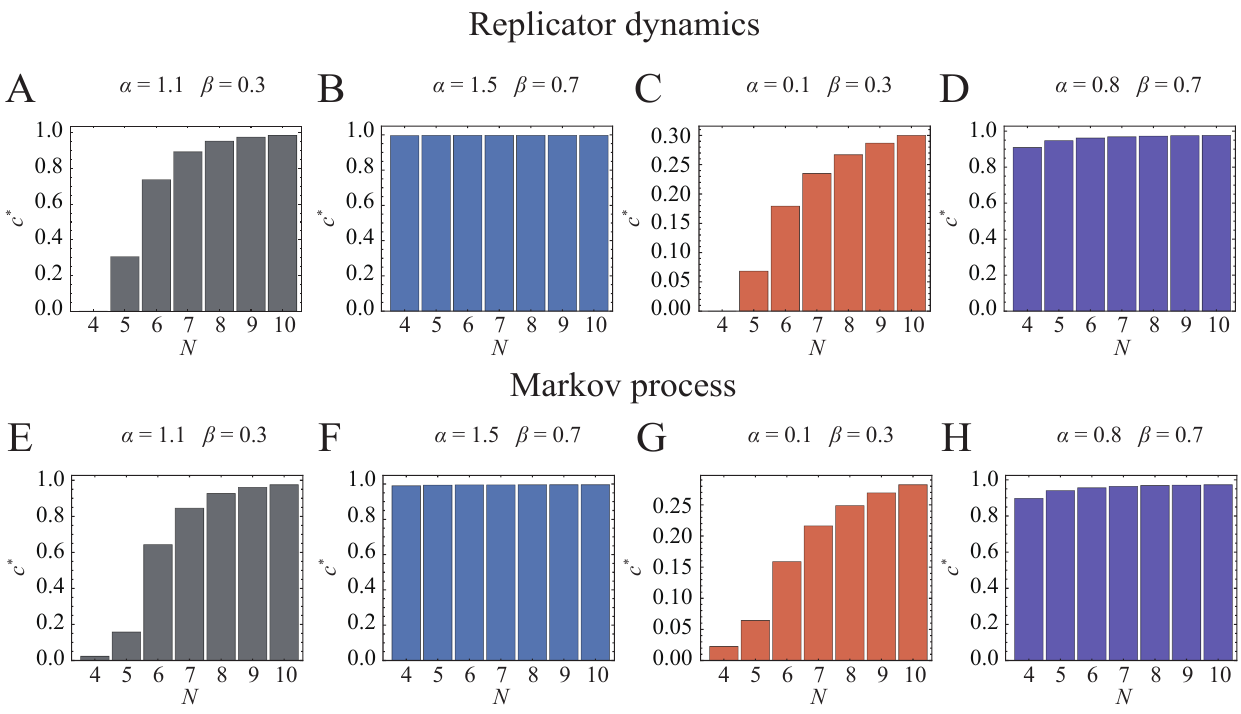}
    \caption{The influence of group size $N$ on the evolutionary dynamics of cooperation under different guidance conditions. For both persuasion-prone ($\alpha>1$) and persuasion-averse ($0<\alpha<1$) populations, regardless of the quality of the guidance environment, the system's cooperation rate increases as group size $N$ grows. Most notably, the enhancement effect of group size expansion on cooperation levels is most pronounced in cases with strong persuasion-prone characteristics and poor guidance conditions. The remaining parameters for panels A-D are $c=1$, $r=n/2$, $\theta=0.01$, and $\delta=0.01$; $Z=100$, $c=1$, $r=n/2$, $\theta=0.01$, $\mu=1/Z$, $\delta=0.01$, and $\gamma=10$ in panels E-H.}\label{n}
\end{figure*}

We further systematically investigated the impact of group size $N$ on the evolution of cooperation  (see Figure \ref{n}). Under the fixed dilemma strength parameter $r=N/2$, both modeling methods reached consistent conclusions: the expansion of group size promotes the emergence of cooperative behavior, a finding that aligns with the social phenomenon of ``many hands make light work." Specifically, we find that in populations with persuasion-prone tendencies ($\alpha>1$) and poor guidance conditions, increasing group size significantly enhances guidance success rates by raising the number of potential guides, thereby driving growth in cooperation rates (see Figure \ref{n} A and E). Besides, it is shown that for populations with persuasion-prone tendencies in favorable guidance environments, relatively high population cooperation rates can be maintained spontaneously even at small group sizes (see Figure \ref{n} B and F). Figure \ref{n} C and G present a more complex scenario: in populations with persuasion-averse tendencies ($\alpha<1$) and poor guidance conditions, although expanding group size can still improve cooperation rates, the degree of improvement is significantly limited. This is primarily because increased guidance attempts trigger resistance among defectors. In contrast, results from Figure \ref{n} D and H confirm that for persuasion-averse populations in good guidance environments, the cooperation rate remains high, and expanding group size can further promote the formation of even higher cooperation levels.

\section{Conclusion}
Cooperative behavior is widespread in nature, but the potential high payoffs of defection strategies always pose an obstacle to collective cooperation. In this work, altruistic guides promote cooperation by actively persuading defectors to change their strategies, even though this process often involves temporary sacrifices of personal interests. How such self-sacrificing guidance behavior can be selected and maintained in the course of evolution is a central theoretical question for understanding the evolutionary dynamics of cooperative behavior \cite{apicella2019evolution,west2021ten,henrich2021origins}.

In this work, we have revealed the dual role of guidance mechanism in the evolution of cooperative behavior: on one hand, it significantly promotes the emergence and sustained stability of cooperative behavior, while on the other hand, the presence of defectors effectively prevents the occurrence of second-order free-riding. Replicator dynamics outcomes demonstrate that the system exhibits bistable characteristics: when the initial proportion of guides is low, insufficient guidance efficiency leads the system inevitably toward complete defection; whereas when the initial proportion of guides is high, the success probability of guidance increases significantly with the number of guides. Under these conditions, the persistent behavioral intervention by guides enables initially defecting individuals to gradually modify their strategies, ultimately driving the population toward a stable equilibrium with high cooperation levels.

We have identified differential effects of guidance mechanisms through comparative analysis of populations with distinct socio-psychological characteristics. The results demonstrate that for persuasion-prone populations, increasing both the initial guidance success rate and persuasion propensity effectively enhances the evolutionary advantage of guidance strategies. Particularly noteworthy is that stronger persuasion tendencies also expand the basin of attraction for co-existence equilibria between guidance and defection strategies in the dynamical system. However, when the initial guidance success rate is low, replicator dynamics and stochastic dynamics yield divergent predictions: the former shows that cooperative behavior can persist given sufficient initial guides when individuals are susceptible to repeated persuasion, while the latter indicates that even with pronounced persuasion propensity, the ultimate cooperation rate remains difficult to elevate. For persuasion-averse populations, replicator dynamics and stochastic dynamics produce relatively consistent conclusions. We find that reducing aversion to persuasion does facilitate higher population cooperation levels, but the most critical determinant remains the initial guidance success rate. When this success rate is low, stable cooperative behavior proves difficult to establish in the population, even with relatively weak individual-level persuasion aversion.

We have further examined the mechanism by which information costs and guidance costs affect the evolution of cooperation. Consistent with theoretical expectations, increases in information acquisition costs and guidance implementation costs significantly inhibit the spread of guidance strategies within groups. It is worth noting that in an environment where initial guidance is poor, high cost may cause guidance strategies to completely lose their evolutionary advantage, resulting in a situation where no one in the group is willing to take on the responsibility of providing guidance. This finding has important policy implications: for governments or other regulatory agencies, providing low-cost or free information service platforms to reduce barriers to information acquisition for agents may be an effective way to maintain the group's spontaneous guidance mechanism. Such intervention measures not only promote the natural emergence of guiding behaviors, but also bring greater long-term benefits to managers by improving the overall level of cooperation, thereby achieving a double improvement in social benefits and governance effectiveness.

We have also investigated how group size $N$ influences the evolutionary dynamics of the system. We have found that as group size increases, the system can sustain a greater number of guides, which enhances cooperation levels in both types of psychologically distinct groups. However, comparative analysis reveals a crucial distinction: in persuasion-prone groups, the increased number of guides significantly improves the success rate of guidance behavior, leading to rapid growth in cooperation levels. In contrast, for persuasion-averse groups, where the guidance success rate declines with repeated attempts, the positive effect of group size expansion on cooperation levels remains relatively limited.

It should be particularly noted that while our work focuses on the cooperative evolutionary patterns of two typical psychological characteristic populations, real-world populations often exhibit more complex distributions of psychological traits. Existing research indicates that psychological heterogeneity among individuals may significantly influence the evolution of cooperation, a direction that warrants further investigation \cite{chen2021effects,santos2021complexity}. Additionally, for persistent defectors, combining guidance mechanisms with punitive measures may serve as an effective supplementary strategy.
\section*{Declaration of Competing Interest}
The authors declare that they have no known competing financial interests or personal relationships that could have appeared to influence the work reported in this paper.

\section*{Data availability}
No data was used for the research described in the article.

\section*{Acknowledgments}
This work was funded by the National Natural Science Foundation of China (Nos. 62306243, and 62406255),  China Postdoctoral Science Foundation (Certificate Number: 2024M762633), and the Humanities and Social Sciences Research Planning Fund of the Ministry of Education (No. 24XJC630006).

\bibliographystyle{elsarticle-num} 
\bibliography{sn-bibliography}


\end{document}